\documentclass[preprint,11pt]{elsarticle}

\newtheorem{theorem}{Theorem}[section]

\newtheorem{lemma}{Lemma}[section]

\usepackage{amssymb}
\usepackage{amsmath}
\usepackage{color}
\begin{document}

\begin{frontmatter}

\title{On Wasserstein-1 distance in the central limit theorem for elephant random walk}
\author{Xiaohui Ma$^a$,\ \ \ Mohamed El Machkouri$^b$, \ \ \  Xiequan Fan$^{a,}$$^*$  }
 \cortext[cor1]{\noindent Corresponding author. \\
  \mbox{ \ \ \ } \textit{E-mail}: fanxiequan@hotmail.com (X. Fan). }
\address[1]{Center for Applied Mathematics, Tianjin University, 300072 Tianjin,  China}
\address[2]{Laboratoire de Math\'ematiques Rapha\"el Salem,
UMR CNRS 6085, Universit\'e de Rouen, France.}

\begin{abstract}
Recently, the elephant random walk has attracted a lot of attentions.  A wide range of literature is available for the asymptotic behavior of the process, such as the central limit theorems,  functional limit theorems and the law of iterated logarithm. However, there is not result concerning Wasserstein-1 distance for the normal approximations.
In this paper, we  show that the Wasserstein-1 distance in the central limit theorem is totally different when   a
memory parameter $p$  belongs to one of the three cases  $0< p  < 1/2,$  $1/2<  p<3/4$ and  $p=3/4.$
\end{abstract}
\begin{keyword}
Elephant random walk; Central limit theorem; The Wasserstein distance
\vspace{0.3cm}
\MSC Primary 60G42; 60F05; 60G50

\end{keyword}

\end{frontmatter}

\section{Introduction and main results}
Random walks, are widely used in theoretical physics to describe phenomena of traveling motion and mass transport.
As a kind of interesting random walk, the elephant random walk   was introduced by Sch{\"u}tz and Trimper \cite{schutz2004elephants} in 2004 in order to study the memory effects in the non-Markovian processes.
The name being inspired by the fact that elephants can remember where they have been.
It is a one-dimensional discrete-time random walk on integers, which has a complete memory of its whole history. The elephant random walk model can be described  as follows.
 It starts at time zero on $0$. At time $n=1$, the elephant $\eta_1$ has a Rademacher $\mathcal{R}(q)$ distribution, where $q$ lies between $[0,1].$
 For $n\geq1$, we define
$$
\eta_{n+1}=\alpha_n \eta_{\beta_n},
$$
where $\alpha_n$ has a Rademacher $\mathcal{R}(p)$ distribution, $p$ lies between $[0,1]$ and $\beta_n$ is uniformly distributed over $\{1, 2,\ldots, n\}$.
Moreover, $\alpha_n$ is independent of $\eta_{\beta_n}$.
The position at time $n+1$ is
$$
S_{n+1}=\sum_{i=1}^{n+1}\eta_i.
$$
In particular, when $p=0$ or $ 1/2$, the elephant random walk reduces to  the classical symmetric random
walk. Since the latter random walk has been well studied, we skip the cases $p=0$ or $ 1/2$.

The elephant random walk  model has drawn a lot of attention in the past few years.
One of the interesting questions concerns the influence of the memory effect   on the long time behavior.
Quite recently, via a method of connection to P\'olya-type urns, Baur and Bertoin \cite{baur2016elephant} derived the functional limit theorems and showed that the limiting process of elephant random walk turns out to be Gaussian  in the diffusive regime  $(0<p<3/4)$  and the critical regime  $(p=3/4)$. But, in  the superdiffusive regime  $(3/4<p\leqslant 1$),   the limiting process of elephant random walk is non-Gaussian.
See also Coletti, Gava and Sch{\"u}tz
\cite{coletti2017central,CGS17} for central limit theorem (CLT) and strong invariance principle of elephant random walk.
With a martingale approach, Bercu \cite{bercu2017martingale} recovered  the  CLT for the elephant random walk as $0\leq p \leq 3/4$.
For $3/4<p\leq 1$,  he also showed that  the elephant random walk converges to a non-normal random variable.
The strong law of large numbers and the law of iterated logarithm  are also discussed in \cite{bercu2017martingale}.
For the extensions of elephant random walk, we refer to 
Bercu and Laulin \cite{bercu2019on},
where Bercu and Laulin have established the CLT for the multi-dimensional elephant random walk.
Though the asymptotic behavior of the elephant random walk, such as  CLT,  functional limit theorems and the law of iterated logarithm,  has been well studied,  there is no result concerning Wasserstein-1 distance for the normal approximations.  Set $a_1=1$ and denote, for $n\geq 2,$
	$$a_n=\frac{\Gamma(n)\Gamma(2p)}{\Gamma(n+2p-1)}    \ \ \ \  \textrm{and} \ \ \ \ v_n=\sum_{i=1}^n a_i^2,$$
	where $\Gamma(x)$ is the  Gamma function.
	Notice that the exact values of $a_n$ and $v_n$ can be easily calculated via computer. In this paper, our main interest is   Wasserstein-1 distance for  CLT of elephant random walk.
Recall that  the Wasserstein-1 distance for    CLT  is defined as follows: $$\|F_{X_n}-\Phi\|_1=\int_{-\infty}^{+\infty}|F_{X_n}(x)-\Phi(x)|\mathrm{d}x,$$
where $F_{X_n}$ and $\Phi$ are the distribution function of the random variable $X_n$ and the standard normal random variable, respectively. The main result of this paper is the following theorem which gives some Wasserstein-1 distance results for the  CLT of elephant random walk. In particular, we get three different convergence rates depending on the value of the memory parameter $p$.

\begin{theorem}\label{thm3.1} Let $n\geq2.$ The following inequalities hold.
\begin{description}
  \item[\textbf{[i]}] If $0 < p < 1/2,$ then
		\begin{equation}\label{p}
		\|F_{ a_nS_n /\sqrt{v_n}}-\Phi\|_1\leq C\frac{\log n}{\sqrt{n}}.\nonumber
		\end{equation}

\item[\textbf{[ii]}] If $1/2<  p<3/4,$ then
		\begin{equation}
		\|F_{ a_nS_n /\sqrt{v_n}}-\Phi\|_1\leq C_p\frac{\log n}{\sqrt{n^{3-4p}}}.\nonumber
		\end{equation}

\item[\textbf{[iii]}] If $p=3/4,$ then
		\begin{equation}
	\|F_{ a_nS_n /\sqrt{v_n}}-\Phi\|_1\leq C\frac{\log \log n}{\sqrt{\log n}}.\nonumber
		\end{equation}	
\end{description}
Here, $C$ is an absolute constant and $C_p$  depends only on $p$.
\end{theorem}

 The main idea in the proof of Theorem \ref{thm3.1} is to introduce a martingale associated to the elephant random walk and to take advantage of powerfull limit theorems from martingale theory. In particular, we are going to apply the CLT established in Lemma \ref{thm2} below.

\section{Proof of Theorem \ref{thm3.1}}\label{seced}
Define the filtration $\mathcal{F}_n=\sigma(\eta_i,\ 1\leq i\leq n)$ and
\begin{equation}\label{fsdfsdf}
M_n=a_nS_n-2q+1.
\end{equation}
Then  $(M_n, \mathcal{F}_n)_{n\geq 1}$   is a martingale. Indeed, it is easy to see that for
\begin{eqnarray}
\mathbf{E}[M_{n+1}| \mathcal{F}_n]&=&\mathbf{E}[a_{n+1}(S_n+\alpha_n\eta_{\beta_n}) -2q+1 |\mathcal{F}_n]\nonumber\\
&=&a_{n+1}S_n+a_{n+1}\mathbf{E}[\alpha_n]\mathbf{E}[\eta_{\beta_n}|\mathcal{F}_n]-2q+1.\nonumber
\end{eqnarray}
Taking into account that
\begin{eqnarray}\label{dssdg}
\mathbf{E}[\alpha_n]=p+ (-1)(1-p) \ \ \  \textrm{ and} \ \ \   \mathbf{E}[\eta_{\beta_n}|\mathcal{F}_n]= \sum_{i=1}^n \frac1n \eta_i,
\end{eqnarray}
we deduce that
\begin{eqnarray}
\mathbf{E}[M_{n+1}| \mathcal{F}_n]
&=&a_{n+1} \Big(S_n  + \big(p+ (-1)(1-p)\big) \sum_{i=1}^n \frac1n \eta_i  \Big)  -2q+1 \nonumber\\
&=&a_{n+1} \Big(S_n  + \big(2p -1 \big)   \frac{S_n}{n} \Big) -2q+1 \nonumber\\
&=&a_nS_n-2q+1\nonumber\\
&=&M_n. \qquad
\end{eqnarray}
Let $( \Delta M_n)_{n\geq1}$ be the martingale differences defined by $ \Delta M_1=M_1$ and, for $n\geq 2,$
$$ \Delta M_n=M_n-M_{n-1}.$$
Denote $\langle M \rangle_n $ be the quadratic variation of $M_n,$ that is $\langle M\rangle_n=\sum_{i=1}^n\mathbf{E}[\Delta M_i^2|\mathcal{F}_{i-1}]$.
The normalized martingale $(M_n/\sqrt{v_n}, \mathcal{F}_n)_{n\geq 1}$  satisfies the following property.
\begin{lemma}\label{lem3.2}
For   $n\geq 2$, we have
	 $$\bigg\|  \frac{\Delta M_i }{\sqrt{v_n}} \bigg\|_\infty \leq   \frac{2}{\sqrt{v_n} }  a_i  $$
and for any $t\geq 1,$  there exists a positive constant $C_t$ such that
\begin{eqnarray*}\displaystyle
 \bigg\|\frac{\langle M \rangle_n}{v_n}-1\bigg\|_t    \leq  \left\{ \begin{array}{ll}\displaystyle
\frac{ C_t }{3-4p}     n^{-1},\ \   & \textrm{if $0< p <3/4  $,}\\
\\
C_t  \big( \log n    \big)^{-1}, \ \   & \textrm{if $  p=3/4$.}
\end{array} \right.
\end{eqnarray*}
\end{lemma}

\noindent\textbf{Proof.}
Clearly,   from Lemma \ref{lem1} of Fan et al.\ \cite{FHM20}, we have  $$\|\Delta M_i/\sqrt{v_n} \|_\infty \leq 2 a_i/\sqrt{v_n},$$
which gives the first desired inequality.
Next we give an estimation of $ \| \langle M\rangle_n/ v_n -1 \|_t,  1\leq t < \infty,$ for   $0<p \leq 3/4.$
From (\ref{fsdfsdf}), we get $\Delta M_k = a_k\varepsilon_k=a_k(S_k-\gamma_{k-1}S_{k-1}). $ Thus, it holds
\begin{eqnarray*}
	\mathbf{E}[(\Delta M_k)^2   |\mathcal{F}_{k-1}]&=&a_k^2 \mathbf{E}[ (S_k-\gamma_{k-1}S_{k-1})^2 |  \mathcal{F}_{k-1}]\\
&=&a_k^2 \big( \mathbf{E}[S_k^2|  \mathcal{F}_{k-1}]  -2\gamma_{k-1}S_{k-1}\mathbf{E}[ S_k |  \mathcal{F}_{k-1}]+\gamma_{k-1}^2S_{k-1}^2 \big) .
\end{eqnarray*}
By (\ref{dssdg}) and $(\alpha_kX_{\beta_k})^2 =1$, it is easy to see that
\begin{eqnarray*}
  \mathbf{E}[S_k^2|  \mathcal{F}_{k-1}] &=& \mathbf{E}[ (S_{k-1}+\alpha_kX_{\beta_k})^2|\mathcal{F}_{k-1}] \\
  &=& S_{k-1}^2+ 2 S_{k-1} \mathbf{E}[  \alpha_k X_{\beta_k} |\mathcal{F}_{k-1}] +1 \\
  &=&S_{k-1}^2+ 2 \frac{2p-1}{k-1}  S_{k-1}^2   +1 \\
  &=&  (2 \gamma_{k-1} -1) S_{k-1}^2   +1
\end{eqnarray*}
and that
\begin{eqnarray*}
\mathbf{E}[ S_k |  \mathcal{F}_{k-1}] &=& \mathbf{E}[  S_{k-1}+\alpha_kX_{\beta_k} |\mathcal{F}_{k-1}]
   =  S_{k-1} +   \frac{2p-1}{k-1}  S_{k-1} =  \gamma_{k-1}   S_{k-1}.
\end{eqnarray*}
Thus, we have $\mathbf{E}[(\Delta M_1)^2   ] =1=a_1^2 $ and for $k\geq2,$
\begin{eqnarray*}
	\mathbf{E}[(\Delta M_k)^2   |\mathcal{F}_{k-1}]
&=&a_k^2 \big( (2 \gamma_{k-1} -1) S_{k-1}^2   +1   -2\gamma_{k-1}^2S_{k-1}^2 +\gamma_{k-1}^2S_{k-1}^2 \big) \\
&=&a_k^2 \big( 1-( \gamma_{k-1} -1)^2 S_{k-1}^2 ) \\
&=& a_k^2   - ( 2p -1)^2 a_k^2 ( \frac{S_{k-1}}{k-1} )^2   .
\end{eqnarray*}
Hence, by the definition of $v_n$ and $M_k$, we obtain
\begin{eqnarray*}
	\langle M\rangle_n&=&v_n-(2p-1)^2\bigg(\sum_{k=1}^{n-1}\Big (\frac{a_{k+1}}{a_k}\Big)^2\Big(\frac{M_k}{k}\Big)^2\bigg).
\end{eqnarray*}
Since $\frac{a_{n+1}}{a_n}\sim 1$ as $n\rightarrow \infty$ (cf.\ equality (\ref{a10})), by Minkowski's inequality,   we have for $t \geq 1,$
\begin{equation}
\|\langle M\rangle_n-v_n  \|_t \leq C (2p-1)^2\Big\|\sum_{k=1}^{n-1}(\frac{M_k}{k})^2\Big\|_t \leq C \sum_{k=1}^{n-1}\frac{1}{k^2}\|M_k\|_{2t}^2,
\end{equation}
where $C$ is a positive constant which values may change from line to line.  Using  Rio's inequality (cf. Theorem 2.1 of \cite{Rio09}) and $\| \Delta M_i  \|_\infty \leq 2 a_i $,   we derive that for $t \geq  1,$
\begin{eqnarray*}
	\|M_k \|_{2t}^2\leq  (2t-1) \sum_{i=1}^k \| \Delta M_i\|_{2t}^2 \leq (2t-1)4   v_k.
\end{eqnarray*}
In the diffusive regime $0< p<3/4 $, by (\ref{a15}), we get for $t\geq 1,$
\begin{equation}
\|\langle M\rangle_n-v_n  \|_t \leq 4C (2t-1)\sum_{k=1}^{n-1}\frac{1}{k^2}v_k \leq C(2t-1 ) \frac{\Gamma{(2p)}^2}{3-4p}\sum_{k=1}^{n-1}k^{1-4p}\leq C \frac{2t-1 }{3-4p}  n^{2-4p}.
\end{equation}
In the critical regime $p=3/4,$ by (\ref{a16}), we have for $t\geq 1,$
\begin{equation}
\|\langle M\rangle_n-v_n  \|_t \leq C(2t-1 )\sum_{k=1}^{n-1}\frac{\log k}{k^2} \leq  C (2t-1) .
\end{equation}
Consequently, again by (\ref{a15}) and (\ref{a16}), we obtain the desired inequality.
This completes the proof of lemma  \ref{lem3.2}.
\hfill\qed\\

For simplicity of notation, denote by
$$\xi_i=\frac{\Delta M_i }{\sqrt{v_n}} , \ \ \  1\leq i \leq  n,\ \ \   X_n=\sum_{i=1}^n \xi_i \ \ \  \textrm{and} \ \ \  \langle X\rangle_n=\sum_{i=1}^n\mathbf{E}[\xi_i^2|\mathcal{F}_{i-1}]. $$
From the last lemma, we have:
\begin{lemma}\label{lem3.3}
For   $n\geq 2$, we have
\begin{eqnarray*}\displaystyle
 \big\| \xi_i \big\|_\infty  \leq     \left\{ \begin{array}{ll}\displaystyle
 C \ n^{-1/2},\ \   & \textrm{if  $0<  p <1/2  $,}\\
C_p \ n^{- (3-4p)/2},\ \   & \textrm{if  $1/2 < p <3/4  $,} \\
\displaystyle    C( \log n)^{-1/2}, \ \   & \textrm{if  $  p=3/4$,}
\end{array} \right.
\end{eqnarray*}
and for any $t\geq 1,$  there exists a positive constant $C_t$ such that
\begin{eqnarray*}\displaystyle
 \Big\|  \langle X\rangle_n -1\Big\|_t    \leq  \left\{ \begin{array}{ll}\displaystyle
\frac{ C_t  }{3-4p}     n^{-1},\ \   & \textrm{if  $0< p <3/4$,}\\
\\
C_t   \big( \log n    \big)^{-1}, \ \   & \textrm{if  $  p=3/4$,}
\end{array} \right.
\end{eqnarray*}
where $C$ is an absolute constant and $C_p$ depends only on $p$.
\end{lemma}

\noindent\textbf{Proof.} We only need to prove the first inequality. The second one holds obviously.
By Stirling's formula $$\log  \Gamma(x)  =(x-\frac12) \log x -x + \frac12 \log 2 \pi + O(\frac{1}{x})  \ \ \ \  \ \textrm{as}\ \ x \rightarrow \infty,$$ we deduce that
\begin{eqnarray}\label{a10}
\lim\limits_{n\to \infty} a_n n^{ 2p-1} =\Gamma(2p).
\end{eqnarray}
Moreover, in the diffusive regime $(0\leq p< 3/4)$, we have
\begin{eqnarray}\label{a15}
\lim\limits_{n\to \infty} \frac{ v_n}{ n^{3-4p}} = \frac{\Gamma{(2p)}^2}{3-4p},
\end{eqnarray}
and, in the critical regime $(p=3/4)$,  it holds
\begin{eqnarray}\label{a16}
\lim\limits_{n\to \infty} \frac{ v_n}{  \log n}    =\frac{3}{4}.
\end{eqnarray}
See also Bercu \cite{bercu2017martingale}  for the equalities (\ref{a10})-(\ref{a16}).
Hence, by Lemma \ref{lem3.2} and the  equalities  (\ref{a10})-(\ref{a16}),  we have
\begin{eqnarray*}\displaystyle
 \big\| \xi_i \big\|_\infty  \leq     \left\{ \begin{array}{ll}\displaystyle
 C \ n^{-1/2},\ \   & \textrm{if $0< p < 1/2  $,}\\
C_p \ n^{- (3-4p)/2},\ \   & \textrm{if $1/2 < p <3/4  $,} \\
\displaystyle    C( \log n)^{-1/2}, \ \   & \textrm{if $  p=3/4$,}
\end{array} \right.
\end{eqnarray*}
which gives the desired inequality.
\hfill\qed

In the sequel   we shall use the following conditions:
\begin{description}
  \item[(A1)]   There exist a positive constant $\rho$ and a number $\varepsilon_n\in (0, \frac{1}{2}]$, such that
$$\mathbf{E}[|\xi_i|^{2+\rho}| \mathcal{F}_{i-1}]\leq \varepsilon_n^{\rho}\mathbf{E}[\xi_i^2| \mathcal{F}_{i-1}]  \quad a.s.\  \mbox{ for all  }i=1, 2,\ldots, n.$$

 \item[(A2)]  There exists a number $\delta_n\in [0,\frac{1}{2}]$, such that $$\|\langle X\rangle_n -1 \|_\infty \leq \delta_n^2 \quad   \mbox{ for all  }i=1, 2,\ldots, n.$$

\end{description}

With conditions (A1) and (A2), we have the following lemma for the mean martingale CLT, which is of independent interest.
The proof of the lemma  is given in the next section.  Notice that the convergence rate of the Wasserstein-1 distance    in the following lemma
is same to the best possible Berry-Esseen bound for martingales, see  Fan \cite{fan2019exact} and El Machkouri and Ouchti \cite{el2007exact}.

\begin{lemma}\label{thm1}
 Assume the conditions (A1) and (A2). Then
\begin{equation}\label{sgfdf}
\|F_{X_n}-\Phi\|_1\leq C_{\rho}(\hat{\varepsilon}_n+\delta_n),
\end{equation}
where
\begin{eqnarray}\label{thma}
\hat{\varepsilon}_n  = \left\{ \begin{array}{l}
{\varepsilon_n ^\rho } \ \ \  \qquad \ \ \ \ \mbox{ if }  \rho \in (0, 1),\\
\varepsilon_n |\log \varepsilon_n| \ \ \ \ \mbox{ if } \rho  \geq 1.
\end{array} \right.
\end{eqnarray}
\end{lemma}

The condition (A2) is very restrictive, while we have the following lemma when  condition (A2)  is dropped. The proof of this lemma is present in Section \ref{dgf}.
\begin{lemma}\label{thm2}
  Assume condition (A1). Then for any $t\geq 1$,
\begin{equation}\label{fsx}
  \|F_{X_n}-\Phi\|_1\leq C_{t,\rho}\Big(\hat{\varepsilon}_n+\Big(\mathbf{E}|\langle X\rangle_n -1|^t+\mathbf{E}\max_{1\leq i\leq n}|\xi_i|^{2t}\Big)^{1/{2t}}\Big),
\end{equation}
where $\hat{\varepsilon}_n$ is defined by (\ref{thma}).
\end{lemma}

Now, we are in position to prove  Theorem \ref{thm3.1}.
 Recall that $X_n=(a_n S_n -2q+1)/\sqrt{v_n}. $
By  Lemma  \ref{lem3.3}, it is easy to see that
\begin{eqnarray*}\displaystyle
	\mathbf{E}|\langle X\rangle_n -1| \leq  \left\{ \begin{array}{ll}\displaystyle
\frac{ C   }{3-4p}     n^{-1},\ \   & \textrm{if $0< p <3/4$,}\\
\\
C   \big( \log n    \big)^{-1}, \ \   & \textrm{if $  p=3/4$.}
	\end{array} \right.
\end{eqnarray*}
Applying Lemmas \ref{lem3.3} and \ref{thm2} to $X_n$ with $t=1$ and
\begin{eqnarray*}\displaystyle
	\varepsilon_n  =   \left\{ \begin{array}{ll}\displaystyle
		C \ n^{-1/2},\ \   & \textrm{if $0< p < 1/2  $,}\\
		C_p \ n^{- (3-4p)/2},\ \   & \textrm{if $1/2 < p <3/4  $,} \\
		\displaystyle    C( \log n)^{-1/2}, \ \   & \textrm{if $  p=3/4$,}
	\end{array} \right.
\end{eqnarray*}
 we  obtain the desired inequalities  of Theorem \ref{thm3.1}. Note that  \begin{equation*}
 \Vert F_{a_nS_{n}/\sqrt{v_{n}}}-\Phi{\Vert}_{1} \leq \Vert F_{(a_nS_{n}-2q+1)/\sqrt{v_{n}}}-\Phi{\Vert}_{1} +1/\sqrt{v_{n}}.
\end{equation*}
 This completes the proof of Theorem \ref{thm3.1}.

\hfill\qed

\section{Proof  of Lemma \ref{thm1}}\label{s2}

 In the sequel,  constants $C$ and $C_p$ are always numerical constants that may change between appearances.
In the proof  of Lemma \ref{thm1}, we shall  use the following technical lemma  of Van Dung et al.\ \cite{son20141}.
\begin{lemma}\label{lem1}
Let $X$ and $Y$ be two random variables. Then for $r>1/2,$
$$\|F_{X}-\Phi\|_1\leq \|F_{X+Y}-\Phi\|_1+2(2r+1)\|\mathbf{E}[Y^{2r}|X]\|_{\infty}^{1/{2r}}.$$
\end{lemma}

The following   lemma  can be found in Fan \cite{fan2019exact}.
\begin{lemma}\label{lf}
Under condition (A1), $\xi_i$ has a bounded conditional variance, that is,
\begin{equation}\label{fan}
\mathbf{E}[\xi_i^2|\mathcal{F}_{i-1}]\leq \varepsilon_n^2 \quad a.s.
\end{equation}
for all $i=1, 2,\ldots, n.$
\end{lemma}

In the proof  of Lemma \ref{thm1},
we also need the following lemma.
\begin{lemma}\label{lem2}
  Let $G$ be a function $\mathbb{R}\to \mathbb{R}$, which has derivative $G'$ which together with $G$ belong to $L^1(\mathbb{R})$.
  For a random variable $X$ and a constant $a>0$, it holds
  $$\mathbf{E}\bigg[\int_{-\infty}^{+\infty}G\Big(\frac{X+t}{a}\Big)\mathrm{d}t\bigg]\leq \|G'\|_1\|F_{X}-\Phi\|_1+a\|G\|_1.$$
\end{lemma}

\noindent\textbf{Proof.}
Let $N$ be a standard normally distributed random variable. By the definition of expectation, it is easy to see that
\begin{flalign}\label{la}
  &\bigg|\mathbf{E}\bigg[\int_{-\infty}^{+\infty}G\Big(\frac{X+t}{a}\Big)\mathrm{d}t\bigg]-\mathbf{E}\bigg[\int_{-\infty}^{+\infty}G\Big(\frac{N+t}{a}\Big)\mathrm{d}t\bigg]\bigg|\nonumber\\
  &\ \ \ \ \ \ \ \ \ \ \ \ \ \  =\frac{1}{a}\bigg|\int_{-\infty}^{+\infty}\int_{-\infty}^{+\infty}\Big(F_X(x)-\Phi(x)\Big)G'\Big(\frac{x+t}{a}\Big)\mathrm{d}t\mathrm{d}x\bigg|.
\end{flalign}
Using integration by substitution with  $ \left\{ \begin{array}{l}
u = \frac{{x + t}}{a}\\
v = x
\end{array} \right.$,
we deduce that the r.h.s. of equality (\ref{la}) equals to $\big|\int_{-\infty}^{+\infty}\int_{-\infty}^{+\infty}(F_X(v)-\Phi(v))G'(u)\mathrm{d}u\mathrm{d}v\big|.$
Hence, from (\ref{la}), we have
\begin{eqnarray}\label{lb}
 \bigg|\mathbf{E}\bigg[\int_{-\infty}^{+\infty}G\Big(\frac{X+t}{a}\Big)\mathrm{d}t\bigg]-\mathbf{E}\bigg[\int_{-\infty}^{+\infty}G\Big(\frac{N+t}{a}\Big)\mathrm{d}t\bigg]\bigg|
   \leq \|G'\|_1\|F_{X}-\Phi\|_1.
\end{eqnarray}
On the other hand, it holds
\begin{equation}\label{lc}
  \mathbf{E}\bigg[\int_{-\infty}^{+\infty}G\Big(\frac{N+t}{a}\Big)\mathrm{d}t\bigg]=\int_{-\infty}^{+\infty}\int_{-\infty}^{+\infty}G\Big(\frac{x+t}{a}\Big)\Phi'(x)\mathrm{d}t\mathrm{d}x
\leq a\|G\|_1.
\end{equation}
Combining  (\ref{lb}) and (\ref{lc}) together, we obtain the desired inequality.
\hfill\qed \\

Now, we are  in position to prove Lemma \ref{thm1}.
We follow the argument of Grama and Haeusler \cite{grama2000large},
where Grama and Haeusler proved a Berry-Esseen's bound for CLT instead of Wasserstein-1 distance for CLT.
For simplicity of notations, denote
\begin{eqnarray}\label{ga}
\gamma_n  = \left\{ \begin{array}{l}
\varepsilon_n+\delta_n  \qquad \ \ \ \ \mbox{ if }  \rho \in (0, 1),\\
\varepsilon_n |\log \varepsilon_n|+\delta_n\  \ \mbox{if } \rho  \geq 1.
\end{array} \right.
\end{eqnarray}
Set $T=1+\delta_n$, and introduce a modification of the conditional variance $\langle X \rangle$ as follows:
\begin{equation}
V_k=\langle X \rangle_k\mathbf{1}_{\{k<n\}}+T\mathbf{1}_{\{k=n\}}.
\end{equation}
It is obvious that $V_k=\langle X \rangle_k$, $1\leq k<n$, $V_n=T \geq \langle X \rangle_n.$
Let $c$ be a large constant depending on $\rho$ and satisfies $c\geq2$, whose value will be chosen later. We introduce a discrete time predictable process
$$ A_k=c^2\gamma_n^2+T-V_k, \qquad k=1, \ldots, n.$$
Then $A_k$ is $\mathcal{F}_{k-1}$-measurable and non-increasing, which is easily to verify.
Moreover, it is obvious that $A_0=c^2\gamma_n^2+T$, $A_n=c^2\gamma_n^2$.
For $y>0$ and $u,x \in \mathbb{R}$, denote
\begin{equation*}
  \Phi_u(x,y)=\Phi\bigg(\frac{u-x}{\sqrt{y}}\bigg).
\end{equation*}
Let $N$ be a standard normal  random variable, which is independent of $X_n$.
By Lemma \ref{lem1} we deduce that
\begin{flalign}
\|F_{X_n}-\Phi\|_1
&\leq \|F_{X_n+c\gamma_n N}-\Phi\|_1+C_1\gamma_n\nonumber\\
&=\int_{-\infty}^{+\infty}\Big|\mathbf{E}\Phi_u(X_n,A_n)-\Phi(u)\Big|\mathrm{d}u+C_1\gamma_n\nonumber.
\end{flalign}
Applying  triangle inequality to the first term in the last line, we have
\begin{flalign}
\|F_{X_n}-\Phi\|_1
&\leq\int_{-\infty}^{+\infty}\Big|\mathbf{E}\Phi_u(X_n,A_n)-\mathbf{E}\Phi_u(X_0,A_0)\Big|\mathrm{d}u
+\int_{-\infty}^{+\infty}\Big|\mathbf{E}\Phi_u(X_0,A_0)-\Phi(u)\Big|\mathrm{d}u+C_1\gamma_n.\nonumber
\end{flalign}
For the second term in the r.h.s.\ of the last inequality,  by the definition of $A_0,$ it is easy to see that
\begin{equation}
  \int_{-\infty}^{+\infty}\Big|\mathbf{E}\Phi_u(X_0,A_0)-\Phi(u)\Big|\mathrm{d}u=\int_{-\infty}^{+\infty}\Big|\Phi(\frac{u}{\sqrt{c^2\gamma_n^2+T}})-\Phi(u)\Big|\mathrm{d}u\leq C_2\gamma_n,
\end{equation}
then it follows that
\begin{equation}\label{g}
  \|F_{X_n}-\Phi\|_1\leq\int_{-\infty}^{+\infty}\Big|\mathbf{E}\Phi_u(X_n,A_n)-\mathbf{E}\Phi_u(X_0,A_0)\Big|\mathrm{d}u+C_3\gamma_n.
\end{equation}

Next, we give an estimation for   $\int_{-\infty}^{\infty}|\mathbf{E}\Phi_u(X_n,A_n)-\mathbf{E}\Phi_u(X_0,A_0)|\mathrm{d}u.$
By a simple telescope, we derive that
\begin{equation*}
  \mathbf{E}\Phi_u(X_n,A_n)-\mathbf{E}\Phi_u(X_0,A_0)=\mathbf{E}\Big[\sum_{k=1}^n\big(\Phi_u(X_k, A_k)-\Phi_u(X_{k-1}, A_{k-1})\big)\Big].
\end{equation*}
Using the fact
\begin{equation*}
  \frac{\partial^2}{\partial x^2}\Phi_u(x,y)=2\frac{\partial}{\partial y}\Phi_u(x,y),
\end{equation*}
we deduce that
\begin{equation}\label{ab}
   \mathbf{E}\Phi_u(X_n,A_n)-\mathbf{E}\Phi_u(X_0,A_0)=I_1-I_2+I_3,
\end{equation}
where
\begin{eqnarray}
  I_1&=&\mathbf{E}\Big[\sum_{k=1}^n\Phi_u(X_k, A_k)-\Phi_u(X_{k-1}, A_{k})-\frac{\partial}{\partial x}\Phi_u(X_{k-1},A_k)\xi_k-\frac{1}{2}\frac{\partial^2}{\partial x^2}\Phi_u(X_{k-1},A_k)\xi_k^2 \Big],\\
  I_2&=&\mathbf{E}\Big[\sum_{k=1}^n\Phi_u(X_{k-1}, A_{k-1})-\Phi_u(X_{k-1}, A_k)-\frac{\partial}{\partial y}\Phi_u(X_{k-1},A_k)\Delta V_k \Big],\\
  I_3&=&\frac{1}{2}\mathbf{E}\Big[\sum_{k=1}^n\frac{\partial^2}{\partial x^2}\Phi_u(X_{k-1},A_k)(\Delta\langle X\rangle_k-\Delta V_k)\Big],
\end{eqnarray}
with $\Delta\langle X\rangle_k=\langle X\rangle_k-\langle X\rangle_{k-1}$ and $ \Delta V_k=V_k-V_{k-1}, k=1, 2,\ldots, n.$
Notice that, as $\frac{\partial^2}{\partial x^2}\Phi_u(X_{k-1},A_k)$ is $\mathcal{F}_{k-1}$-measurable, it holds
$$\mathbf{E}\Big[ \frac{\partial^2}{\partial x^2}\Phi_u(X_{k-1},A_k) \xi_k^2   \Big]  =  \mathbf{E}\Big[ \frac{\partial^2}{\partial x^2}\Phi_u(X_{k-1},A_k) \Delta\langle X\rangle_k \Big].$$
Next, we give estimates for $\int_{-\infty}^{\infty}|I_1|\mathrm{d}u$, $\int_{-\infty}^{\infty}|I_2|\mathrm{d}u$ and $\int_{-\infty}^{\infty}|I_3|\mathrm{d}u$, respectively. Denote by $\varphi$  the density function of the standard normal random variable. Moreover, $\theta_k$'s stand for  random variables satisfying $0\leq\theta_k\leq1,$ which may take different values in different places.

\textbf{(a)} Estimate of $\int_{-\infty}^{\infty}|I_1|\mathrm{d}u$.

For brevity of notation,  denote by $T_{k-1}=\frac{u-X_{k-1}}{\sqrt{A_k}}$. Clearly, it holds
$$  I_1 = \sum_{k=1}^n  \mathbf{E}  R_k \, ,  $$
where
\begin{eqnarray}\label{r}
 R_{k}&:=& \Phi_u(X_k, A_k)-\Phi_u(X_{k-1}, A_{k})-\frac{\partial}{\partial x}\Phi_u(X_{k-1},A_k)\xi_k-\frac{1}{2}\frac{\partial^2}{\partial x^2}\Phi_u(X_{k-1},A_k)\xi_k^2\nonumber\\
  &=&\Phi(T_{k-1}-\frac{\xi_k}{\sqrt{A_k}})-\Phi(T_{k-1})-\Phi'(T_{k-1})(-\frac{\xi_k}{\sqrt{A_k}})-\frac{1}{2}\Phi''(T_{k-1})(\frac{\xi_k}{\sqrt{A_k}})^2.
\end{eqnarray}
To bound up the bound (\ref{r}), we distinguish two cases as follows.

(1) Case of $|\frac{\xi_k}{\sqrt{A_k}}|\leq 1+\frac{|T_{k-1}|}{2}$.

When $|\frac{\xi_k}{\sqrt{A_k}}|\leq 1 $, by a three term Taylor's expansion,   we have
\begin{eqnarray*}
  |R_k|&=&\Big|\frac{1}{3!}\Phi'''(T_{k-1}-\theta_k\frac{\xi_k}{\sqrt{A_k}})(\frac{\xi_k}{\sqrt{A_k}})^3 \Big|\nonumber\\
  &\leq& \Big|\Phi'''(T_{k-1}-\theta_k\frac{\xi_k}{\sqrt{A_k}})\Big|\Big|\frac{\xi_k}{\sqrt{A_k}}\Big|^{2+\rho}.
\end{eqnarray*}
When $|\frac{\xi_k}{\sqrt{A_k}}|> 1$, by a two term Taylor's expansion,   we get
\begin{eqnarray*}
  |R_k|&\leq&\frac{1}{2}\Big( |\Phi''(T_{k-1}-\widehat{\theta}_k\frac{\xi_k}{\sqrt{A_k}})| +|\Phi''(T_{k-1})|\Big)(\frac{\xi_k}{\sqrt{A_k}})^2\nonumber\\
  & \leq &\Big|\Phi''(T_{k-1}-\theta_k'\frac{\xi_k}{\sqrt{A_k}})\Big| \Big|\frac{\xi_k}{\sqrt{A_k}}\Big|^{2+\rho},
\end{eqnarray*}
where
\begin{eqnarray*}
\theta_k'  = \left\{ \begin{array}{l}
	\widehat{\theta}_k,  \quad  \mbox{ if }|\Phi''(T_{k-1}-\widehat{\theta}_k\frac{\xi_k}{\sqrt{A_k}})|\geq |\Phi''(T_{k-1})|  ,\\
	0,\ \ \ \ \mbox{ otherwise} .
\end{array} \right.
\end{eqnarray*}
Recall that $|\frac{\xi_k}{\sqrt{A_k}}|\leq 1+\frac{|T_{k-1}|}{2}$ and $0\leq\theta_k, \theta_k'\leq1$.
By the inequality $\max\{|\Phi''(t)|, |\Phi'''(t)|\}\leq \varphi(t)(1+t^2)$,    it follows that
\begin{eqnarray}\label{rr}
  |R_k|&\leq&\varphi(T_{k-1}-\theta_1\frac{\xi_k}{\sqrt{A_k}})\Big(1+(T_{k-1}-\theta_1\frac{\xi_k}{\sqrt{A_k}})^2\Big)\Big|\frac{\xi_k}{\sqrt{A_k}}\Big|^{2+\rho}\nonumber\\
&\leq& g_1(T_{k-1})\Big|\frac{\xi_k}{\sqrt{A_k}} \Big|^{2+\rho},
\end{eqnarray}
where
\begin{equation*}
  g_1(z)=\sup_{|t-z|\leq 1+|z|/2}\varphi(t)(1+t^2).
\end{equation*}

(2) Case of $|\frac{\xi_k}{\sqrt{A_k}}|> 1+\frac{|T_{k-1}|}{2}$.

It is easy to see that for $|\Delta x|>1+\frac{|x|}{2},$ it holds
\begin{flalign}
  \Big|\Phi(x+\Delta x)-\Phi(x)-\Phi'(x)-\frac{1}{2}\Phi''(x)(\Delta x)^2\Big|
  &\leq\Big (\frac{|\Phi(x+\Delta x)-\Phi(x)|}{|\Delta x|^{2+\rho}}+|\Phi'(x)|+|\Phi''(x)|\Big)|\Delta x|^{2+\rho}\nonumber\\
  &\leq \Big(8\frac{|\Phi(x+\Delta x)-\Phi(x)|}{(2+| x|)^{2+\rho}}+|\Phi'(x)|+|\Phi''(x)|\Big)|\Delta x|^{2+\rho}\nonumber\\
  &\leq \Big(\frac{C}{(4+ x^2)^{1+\rho/2}}+|\Phi'(x)|+|\Phi''(x)|\Big)|\Delta x|^{2+\rho}\nonumber\\
  &\leq \frac{C_4}{(4+ x^2)^{1+\rho/2}}|\Delta x|^{2+\rho}.
\end{flalign}
Then  in this case  we have
\begin{equation}\label{rrr}
  |R_k|\leq g_2(T_{k-1})|\frac{\xi_k}{\sqrt{A_k}}|^{2+\rho},
\end{equation}
where $$g_2(z)=\frac{C_4}{(4+z^2)^{1+\rho/2}}.$$

Now, we return to the estimation of  bound (\ref{r}).
Set $$G(z)=g_1(z)+g_2(z),$$
then it is easy to verify that $G\in L^1$, $G'\in L^1$ and $G(z), z> 0,$ is non-increasing in $z$. By (\ref{rr}) and (\ref{rrr}), it holds
\begin{equation}\label{fdfdsg}
  |R_k|\leq G(T_{k-1}) \Big|\frac{\xi_k}{\sqrt{A_k}} \Big|^{2+\rho}.
\end{equation}
By (\ref{fdfdsg}) and condition (A1), that is  $\mathbf{E}[|\xi_k|^{2+\rho}| \mathcal{F}_{k-1}]\leq \varepsilon_n^{\rho}\Delta\langle X \rangle_k,$
 we derive that
\begin{flalign}
  |I_1|&\leq \sum_{k=1}^n \mathbf{E}  |R_k|  \leq   \mathbf{E}\Big[\sum_{k=1}^n G(T_{k-1})\frac{1}{A_k^{1+\rho/2}}\mathbf{E}[|\xi_k|^{2+\rho}| \mathcal{F}_{k-1}]\Big]\nonumber\\
&\leq \varepsilon_n^{\rho}\mathbf{E}\Big[\sum_{k=1}^n G(T_{k-1})\frac{1}{A_k^{1+\rho/2}}\Delta\langle X\rangle_k\Big]\nonumber\\
&\leq \varepsilon_n^{\rho}\mathbf{E}\Big[\sum_{k=1}^n G(T_{k-1})\frac{1}{A_k^{1+\rho/2}}\Delta V_k\Big]. \label{ddgsas}
\end{flalign}
To bound up (\ref{ddgsas}),  we introduce the stopping time $\tau_t$ as follows: for any $t\in [0, T]$,
\begin{equation}
  \tau_t=\min\{k\leq n: \langle X\rangle_k>t\}, \qquad \mbox{where min } \emptyset=n.
\end{equation}
Denote $(\sigma_k)_{k=1,\ldots, n}$ be the increasing sequence of moments while the stopping time has jumps, then $\Delta V_k=\int_{[\sigma_k,\sigma_{k+1})}\mathrm{d}t$ and $k=\tau_t$ for $t\in[\sigma_k,\sigma_{k+1}).$ This follows from the definition of $\tau_t$.
Therefore, we have
\begin{flalign}
  \sum_{k=1}^n G(T_{k-1})\frac{1}{A_k^{1+\rho/2}}\Delta V_k&=\sum_{k=1}^n\int_{[\sigma_k,\sigma_{k+1})}\frac{1}{A_{\tau_t}^{1+\rho/2}}G(T_{\tau_t-1})\mathrm{d}t\nonumber\\
  &=\int_0^T\frac{1}{A_{\tau_t}^{1+\rho/2}}G(T_{\tau_t-1})\mathrm{d}t.    \label{ddgsas02}
\end{flalign}
Set $a_t=c^2\gamma_n^2+T-t.$ By Lemma \ref{lf}, we have    $$\Delta V_n=T-V_{n-1}=1+\delta_n^2-\langle X\rangle_{n-1}=1+\delta_n^2-\langle X\rangle_{n}+\Delta\langle X\rangle_{n}\leq 2\delta_n^2+\varepsilon_n^2.$$   Recall  that $A_k=c^2\gamma_n^2+T-V_k,$ then we deduce that  $A_{\tau_t}\leq a_t$
and
$$A_{\tau_t}  \geq  c^2\gamma_n^2+T-t-\Delta V_{\tau_t} \geq  c^2\gamma_n^2+T-t- (  2\delta_n^2+\varepsilon_n^2 )  \geq  \frac{1}{2}a_t  .$$
Thus, by (\ref{ddgsas}) and (\ref{ddgsas02}) and the fact $ \frac{1}{2}a_t  \leq A_{\tau_t}\leq a_t$, it holds
\begin{flalign}
  |I_1|&\leq \varepsilon_n^{\rho}\mathbf{E}\bigg[\int_0^T\frac{1}{A_{\tau_t}^{1+\rho/2}}G(T_{\tau_t-1})\mathrm{d}t\bigg]\nonumber\\
  &\leq 2^{1+\rho/2}\varepsilon_n^{\rho} \int_0^T\frac{1}{a_t^{1+\rho/2}}\mathbf{E}\Big[G\Big(\frac {u-X_{\tau_t-1}}{\sqrt{a_t}}\Big)\Big]\mathrm{d}t.
\end{flalign}
Applying  Fubini's theorem to the last inequality, we have
\begin{equation}\label{a}
  \int_{-\infty}^{+\infty}|I_1|\mathrm{d}u\leq  2^{1+\rho/2}\varepsilon_n^{\rho} \int_0^T\frac{1}{a_t^{1+\rho/2}}\mathbf{E}\Big[\int_{-\infty}^{+\infty}G\Big(\frac {u-X_{\tau_t-1}}{\sqrt{a_t}}\Big)\mathrm{d}u\Big]\mathrm{d}t.
\end{equation}
By Lemma \ref{lem2}, we get
\begin{equation}\label{b}
  \mathbf{E}\bigg[\int_{-\infty}^{+\infty}G\Big(\frac {u-X_{\tau_t-1}}{\sqrt{a_t}}\Big)\mathrm{d}u\bigg]\leq C_5\|F_{X_{\tau_t-1}}-\Phi\|_1+C_6\sqrt{a_t}.
\end{equation}
By Lemma \ref{lem1}, we have
\begin{equation}\label{c}
  \|F_{X_{\tau_t-1}}-\Phi\|_1\leq \|F_{X_n}-\Phi\|_1+C_7\|\mathbf{E}[(X_n-X_{\tau_t-1})^2|X_{\tau_t-1}]\|_{\infty}^{1/2}.
\end{equation}
Now we consider $\mathbf{E}[(X_n-X_{\tau_t-1})^2|\mathcal{F}_{\tau_t-1}] .$ It is easy to see that
\begin{flalign}\label{d}
  \mathbf{E}[(X_n-X_{\tau_t-1})^2|\mathcal{F}_{\tau_t-1}]&=\mathbf{E}\Big[\sum_{k=\tau_t}^n\xi_k^2|\mathcal{F}_{\tau_t-1}\Big]
  =\mathbf{E}\Big[\sum_{k=\tau_t}^n\mathbf{E}[\xi_k^2|\mathcal{F}_{k-1}]\Big|\mathcal{F}_{\tau_t-1}\Big]\nonumber\\
  &=\mathbf{E}[\langle X\rangle_n-\langle X\rangle_{\tau_t-1}|\mathcal{F}_{\tau_t-1}] \nonumber\\
  &\leq V_n-V_{\tau_t-1}\nonumber\\
  &\leq a_t,
\end{flalign}
where the last line follows by the fact
\begin{equation*}
V_n-V_{\tau_t-1}\leq V_n-V_{\tau_t}+\Delta V_{\tau_t}\leq T-t+\varepsilon_n^2+2\delta_n^2\leq a_t.
\end{equation*}
Combining (\ref{a})-(\ref{d}) together, we obtain
\begin{equation}
  \int_{-\infty}^{+\infty}|I_1|\mathrm{d}u\leq C_8\,\varepsilon_n^{\rho}\int_0^T\frac{\mathrm{d}t}{a_t^{1+\rho/2}}\|F_{X_n}-\Phi\|_1+ C_{9}\,\varepsilon_n ^{\rho} \int_0^T\frac{\mathrm{d}t}{a_t^{(1+\rho)/2}}.
\end{equation}
After some simple calculations we deduce that
\begin{equation*}
  \int_0^T\frac{\mathrm{d}t}{a_t^{1+\rho/2}}\leq \frac{2}{\rho c^{\rho}\gamma_n^{\rho}}
\end{equation*}
and
$$\int_0^T\frac{\mathrm{d}t}{a_t^{(1+\rho)/2}}  \leq \left\{ \begin{array}{l}
C_\rho  \qquad \ \ \ \ \mbox{ if }  \rho \in (0, 1),\\
C |\log \varepsilon_n|\ \ \ \mbox{if } \rho  \geq 1.
\end{array} \right.$$
Since $\varepsilon_n\leq \gamma_n$,
finally we get
\begin{equation}\label{e}
  \int_{-\infty}^{+\infty}|I_1|\mathrm{d}u\leq \frac{C_{10}}{\rho c^\rho}\|F_{X_n}-\Phi\|_1+C_{\rho,1} \hat{\varepsilon}_n ,
\end{equation}
where $\hat{\varepsilon}_n$ is defined in (\ref{thma}) and $\gamma_n$ is defined in (\ref{ga}).

\textbf{(b)} Estimate of $\int_{-\infty}^{+\infty}|I_2|\mathrm{d}u$.

By a two-term Taylor expansion,
 we get
\begin{eqnarray}
  I_2&=&\frac{1}{2}\mathbf{E}\Big[\sum_{k=1}^n\frac{\partial^2}{\partial y^2}\Phi_u(X_{k-1}, A_k-\theta_k\Delta A_k)\Delta A_k^2\Big] \nonumber.
\end{eqnarray}
By  the fact that $\frac{\partial^4}{\partial x^4}\Phi_u(x,y) = 4\frac{\partial^2}{\partial y^2}\Phi_u(x,y),$
 we deduce that
\begin{eqnarray}
  I_2 &=&\frac{1}{8}\mathbf{E}\Big[\sum_{k=1}^n\frac{\partial^4}{\partial x^4}\Phi_u(X_{k-1}, A_k-\theta_k\Delta A_k)\Delta A_k^2\Big]\nonumber \\
     &=&\frac{1}{8 \sqrt{2 \pi}}\mathbf{E}\Big[\sum_{k=1}^n\frac{1}{(A_k-\theta_k\Delta A_k)^2}\varphi'''\Big(\frac{u-X_{k-1}}{\sqrt{A_k-\theta_k\Delta A_k}}\Big)\Delta A_k^2\Big]. \nonumber
\end{eqnarray}
Since $|\Delta A_k|:=A_k-A_{k-1}=\Delta V_k\leq\varepsilon_n^2+2\delta_n^2$, we have
\begin{equation*}
  A_k\leq A_k-\theta_k\Delta A_k\leq c^2\gamma_n^2+T-V_k+\varepsilon_n^2+2\delta_n^2\leq 2A_k.
\end{equation*}
Set $\bar{G}(z)=\sup_{|t-z|\leq 2}|\varphi'''(t)|$, then
$\bar{G}(x z)$ is decreasing in $x> 0$.
It follows that
\begin{eqnarray}
  |I_2|&\leq&  \mathbf{E}\Big[\sum_{k=1}^n\frac{1}{A_k^2}\bar{G}\Big(\frac{u-X_{k-1}}{\sqrt{2A_k}}\Big)\Delta A_k^2\Big] \nonumber \\ \nonumber
  &\leq&  (\varepsilon_n^2+2\delta_n^2)\mathbf{E}\Big[\sum_{k=1}^n\frac{1}{A_k^2}\bar{G}\Big(\frac{T_{k-1}}{\sqrt{2}}\Big)\Delta V_k\Big],
\end{eqnarray}
where the last line follows by the facts that $|\Delta A_k| \leq\varepsilon_n^2+2\delta_n^2$ and  $\Delta A_k=\Delta V_k$.
It  is easy to verify  that
$\bar{G}, \bar{G}'\in L^1$ and $|\varphi'''(z)|\leq \bar{G}(z)$ for any $z\geq 0$.
By an argument similar to the proof of (\ref{e}), we  deduce that
\begin{equation}\label{f}
  \int_{-\infty}^{+\infty}|I_2|\mathrm{d}u\leq \frac{C_{11}}{c}\|F_{X_n}-\Phi\|_1+C_{\rho,2} \hat{\varepsilon}_n .
\end{equation}

\textbf{(c)} Estimate of $\int_{-\infty}^{+\infty}|I_3|\mathrm{d}u$.

Since $\Delta\langle X\rangle_k=\Delta V_k, 1\leq k<n, \Delta V_n-\Delta \langle X\rangle_n\leq 2\delta_n^2,$ we have
\begin{eqnarray*}	
|I_3|&=&\bigg|\frac{1}{2}\mathbf{E}[\frac{\partial^2}{\partial x^2}\Phi_u(X_{n-1},A_n)(\Delta\langle X\rangle_n-\Delta V_n)]\bigg|\\ \nonumber
&\leq& \mathbf{E}[\frac{\delta_n^2}{A_n}|\varphi'(T_{n-1})|]\leq \frac{1}{c^2} \mathbf{E}[|\varphi'(T_{n-1})|].
\end{eqnarray*}
Set $\tilde{G}(z)=\sup_{|z-t|\leq 1}|\varphi'(t)|$, then
\begin{equation*}
|I_3|\leq \frac{1}{c^2}\mathbf{E}[\tilde{G}(T_{n-1})].
\end{equation*}
Notice that $\tilde{G}, \tilde{G}' \in L^1$, by Lemmas \ref{lem2}
and \ref{lem1}, we get
\begin{eqnarray}\label{h}
\int_{-\infty}^{+\infty}|I_3|\mathrm{d}u
&\leq& \frac{1}{c^2}\mathbf{E}\Big[\int_{-\infty}^{+\infty}\tilde{G}\Big(\frac {u-X_{n-1}}{\sqrt{A_n}}\Big)\mathrm{d}u\Big]\nonumber\\
&\leq&\frac{C_{12}}{c^2}\|F_{X_n}-\Phi\|_1+C_{\rho,3} (\hat{\varepsilon}_n+\delta_n) ,
\end{eqnarray}
where the last line follows by the fact that $\sqrt{A_n}=c \, \gamma_n= c\, (\hat{\varepsilon}_n+\delta_n)$.

Combining (\ref{e})-(\ref{h}) together, we get
\begin{equation}
\int_{-\infty}^{+\infty}|I_1|+|I_2|+|I_3|\mathrm{d}u\leq \frac{C_{13}}{\rho c^{\rho}}\|F_{X_n}-\Phi\|_1+C_{\rho,4}(\hat{\varepsilon}_n+\delta_n).
\end{equation}
Implementing the last bound in (\ref{g}), by (\ref{ab}) we derive that
\begin{equation}
  \|F_{X_n}-\Phi\|_1\leq \frac{C_{14}}{\rho c^{\rho}}\|F_{X_n}-\Phi\|_1+C_{\rho,5}(\hat{\varepsilon}_n+\delta_n).
\end{equation}
Choosing the value $c$ such that  $\rho c^{\rho}=2C_{14},$ it follows that
\begin{equation}
  \|F_{X_n}-\Phi\|_1\leq 2C_{\rho,5}(\hat{\varepsilon}_n+\delta_n).
\end{equation}
This completes the proof of Lemma \ref{thm1}.\hfill\qed

\section{Proof of Lemma \ref{thm2}}\label{dgf}
  We follow the approach of Bolthausen \cite{bolthausen1982exact}. The main idea is to construct a new martingale difference  sequence $(\hat{\xi}_i, \hat{\mathcal{F}}_i)_{1\leq i\leq N},$   based on a modification of $(\xi_i, \mathcal{F}_i)_{1\leq i\leq n},$ such that $\sum_{i=1}^{N}\mathbf{E}[\hat{\xi}_i^2|\hat{\mathcal{F}}_{i-1}]=1$ a.s.,
 and then apply Lemma  \ref{lem1}   to the new martingale different sequence for obtaining Wasserstein-1 distance in the central limit theorem for $X_n$, that is
\begin{flalign}
\|F_{X_n}-\Phi\|_1
& \leq \|F_{\hat{X}_N}-\Phi\|_1+2(2t+1)\big\|\mathbf{E}\big[|\hat{X}_N-X_n|^{2t}\big|X_n\big]\big\|_{\infty}^{1/{2t}}.
\end{flalign}
The first term can be estimated via Lemma  \ref{thm1}, and the second term can be dominated via Theorem 2.11 of Hall  and Heyde \cite{HH80}.

We introduce the following stopping time
\begin{equation*}
\tau=\sup\{k\leq n: \langle X\rangle_k\leq 1\}.
\end{equation*}
Let $\delta $ be a positive number  such that $0<\delta\leq\varepsilon_n.$ Let $r=\lfloor\frac{1-\langle X\rangle_\tau}{\delta^2}\rfloor$, where $\lfloor x\rfloor$ stands for the largest integer not exceeding $x.$ Clearly, it holds $r\leq\lfloor \frac{1}{\delta^2}\rfloor.$
Let $N=n+r+1.$ For $\tau+1\leq i\leq \tau +r$,   let $\zeta_i$ be   random variables such that
\begin{equation*}
\mathbf{P}(\zeta_i=\pm \delta \   | \ \mathcal{F}_{\tau})=\frac{1}{2};
\end{equation*}
and  for $i=\tau+r+1$,  let $\zeta_{\tau+r+1}$ be such that
\begin{equation*}
\mathbf{P}\Big(\zeta_{\tau+r+1}=\pm (1-\langle X\rangle_\tau -r\delta^2)^{1/2}  \  \Big| \ \mathcal{F}_{\tau} \Big)= \frac{1}{2 } ,
\end{equation*}
 with the sign determined independent of everything else; and let $\zeta_i=0$ if $\tau+r+1<i\leq N$.
Denote
\begin{equation}
\hat{\xi}_i=\xi_i\mathbf{1}_{\{i\leq \tau\}}+\zeta_i\mathbf{1}_{\{\tau<i\leq \tau+r\}}+\zeta_i\mathbf{1}_{\{i=\tau+r+1\}}, \qquad \  i=1,\ldots, N,
\end{equation}
$\hat{\mathcal{F}}_i=\mathcal{F}_i $ for $i\leq \tau $ and $\hat{\mathcal{F}}_i=\sigma\{\mathcal{F}_\tau, \zeta_j, \tau+1\leq j\leq i\}$ for $\tau+1\leq i\leq N$. 
Then $(\hat{\xi}_i, \hat{\mathcal{F}}_{i})_{1\leq i\leq N}$ still forms a martingale difference sequence.
Moreover, it holds that
\begin{flalign*}
\sum_{i=1}^N\mathbf{E}[\hat{\xi}_i^2|\hat{\mathcal{F}}_{i-1}]
&=\sum_{i=1}^{\tau}\mathbf{E}[\xi_{i}^2|\mathcal{F}_{i-1}]+\sum_{i=\tau+1}^{\tau+r+1}\delta^2+(1-\langle X\rangle_{\tau}-r\delta^2)\nonumber\\
&=1 \qquad a.s.
\end{flalign*}
and that,  by the fact   $\delta\leq \varepsilon_n,$
\begin{equation*}
\mathbf{E}[|\hat{\xi}_i|^{2+\rho}|\hat{\mathcal{F}}_{i-1}]\leq \varepsilon_n^{\rho}\mathbf{E}[\hat{\xi}_i^2|\hat{\mathcal{F}}_{i-1}], \qquad \mbox{for } i=1,\ldots, N.
\end{equation*}
Let $\hat{X}_N=\sum_{i=1}^{N}\hat{\xi}_i,$
by Lemma \ref{thm1}, we have
\begin{equation}
\|F_{\hat{X}_N}-\Phi\|_1\leq C_{\rho}\hat{\varepsilon}_n.
\end{equation}
Using Lemma \ref{lem1}, we obtain for $t> 1/2,$
\begin{flalign}\label{bb}
\|F_{X_n}-\Phi\|_1
& \leq \|F_{\hat{X}_N}-\Phi\|_1+2(2t+1)\Big\|\mathbf{E}\Big[|\hat{X}_N-X_n|^{2t}\Big|X_n\Big]\Big\|_{\infty}^{1/{2t}}\nonumber\\
& \leq C_{\rho}\hat{\varepsilon}_n+C_t\Big(\mathbf{E}|\hat{X}_N-X_n|^{2t}\Big)^{1/{2t}}.
\end{flalign}
For the estimation of $\mathbf{E}|\hat{X}_N-X_n|^{2t},$
we first note that
\begin{equation*}
\hat{X}_N-X_n=\sum_{i=\tau+1}^N(\hat{\xi}_i-\xi_i),
\end{equation*}
where  we put $\xi_i=0 \mbox{ for } i>n.$
As $\tau$ is a stopping time, $(\hat{\xi}_i-\xi_i)_{\tau+1\leq i\leq N}$ still forms a martingale difference sequence. By Theorem 2.11 of Hall  and Heyde \cite{HH80}, it holds for $t>1/2,$
\begin{equation}
\mathbf{E}|\hat{X}_N-X_n|^{2t}\leq C_t\bigg(\mathbf{E}\Big|\sum_{i=\tau+1}^N\mathbf{E}[(\hat{\xi}_i-\xi_i)^2|\hat{\mathcal{F}}_{i-1}]\Big|^t+\mathbf{E}[\max_{\tau+1\leq i\leq N}|\hat{\xi}_i-\xi_i|^{2t}]\bigg). \label{ineq48}
\end{equation}
As $\mathbf{E}[\xi_i\hat{\xi}_i|\hat{\mathcal{F}}_{i-1}]=0$ for $\tau+1\leq i\leq N,$ we get
\begin{flalign*}
\sum_{i=\tau+1}^N\mathbf{E}[(\hat{\xi}_i-\xi_i)^2|\hat{\mathcal{F}}_{i-1}]=\sum_{i=\tau+1}^N\mathbf{E}[\hat{\xi}_i^2|\hat{\mathcal{F}}_{i-1}]+\sum_{i=\tau+1}^N\mathbf{E}[\xi_i^2|\hat{\mathcal{F}}_{i-1}]=1-2\langle X\rangle_\tau+\langle X\rangle_n.
\end{flalign*}
Notice that $1-\mathbf{E}[\xi_{\tau+1}^2|\mathcal{F}_{\tau}]<\langle X\rangle_\tau.$ Hence, by the inequality
\begin{equation}\label{fsfdsd}
 |a+b|^{t} \leq \max\{2^{t-1}, 1\}\left(|a|^{t}+|b|^{t}\right),\ \ \ \  t > 1/2,
\end{equation}
 and  Jensen's inequality, we deduce that
\begin{flalign}
\Big|\sum_{i=\tau+1}^N\mathbf{E}\big[(\hat{\xi}_i-\xi_i)^2|\hat{\mathcal{F}}_{i-1}\big]\Big|^t &\leq \Big|\langle X\rangle_n-1+2\mathbf{E}[\xi_{\tau+1}^2|\mathcal{F}_{\tau}]\Big|^t\nonumber\\
&\leq C_t\Big(|\langle X\rangle_n-1|^t+(\mathbf{E}[\xi_{\tau+1}^2|\mathcal{F}_{\tau}])^t\Big)\nonumber\\
&\leq C_t\Big(|\langle X\rangle_n-1|^t+\mathbf{E}[|\xi_{\tau+1}|^{2t} \big|\mathcal{F}_{\tau}]\Big).
\end{flalign}
Taking expectations on both sides of the last inequality, we have
\begin{flalign}
\mathbf{E}\Big|\sum_{i=\tau+1}^N\mathbf{E}\big[(\hat{\xi}_i-\xi_i)^2|\hat{\mathcal{F}}_{i-1}\big]\Big|^t
&\leq C_t \Big(\mathbf{E}|\langle X\rangle_n-1|^t+\mathbf{E}|\xi_{\tau+1}|^{2t}\Big)\nonumber\\
&\leq C_t \Big(\mathbf{E}|\langle X\rangle_n-1|^t+\mathbf{E}[\max_{1\leq i\leq n}|\xi_i|^{2t}]\Big).\label{aa}
\end{flalign}
Similarly, by   inequality (\ref{fsfdsd}),  it holds
\begin{equation}
\mathbf{E}[\max_{\tau+1\leq i\leq N}|\hat{\xi}_i-\xi_i|^{2t}]\leq C_t \mathbf{E}[\max_{1\leq i\leq n}|\xi_i|^{2t}+\delta^{2t}].
\end{equation}
Applying the last inequality and (\ref{aa}) to (\ref{ineq48}), we get
\begin{equation}
\mathbf{E}|\hat{X}_N-X_n|^{2t}\leq C_t\Big(\mathbf{E}|\langle X\rangle_n-1|^{t}+\mathbf{E}[\max_{1\leq i\leq n}|\xi_i|^{2t}]+\delta^{2t}\Big).
\end{equation}
Combining the last inequality with (\ref{bb}) and letting $\delta\to 0$, we obtain the desired inequality.\hfill\qed

\section*{Data availability statement}
Data sharing is not applicable to this article as no new data were created or analyzed in this study.

\section*{Acknowledgements}
The authors deeply indebted to the editor and the anonymous referee for their helpful comments.
This work has been partially supported by the National Natural Science Foundation
of China (Grant Nos.\,11601375 and 11971063). This work was also funded by CY Initiative of Excellence
(grant "Investissements d'Avenir" ANR-16-IDEX-0008),
Project "EcoDep" PSI-AAP2020-0000000013, and by the Labex MME-DII (https://labex-mme-dii.u-cergy.fr/).

\end{document}